\documentclass[12pt]{article}
\usepackage{latexsym}
\usepackage{amsfonts}
\newcommand{\binom}[2]{\left(\begin{array}{c}#1\\#2\end{array}\right)}

\def\const{{\rm const}}

\def\CP{{\bf P}}
\def\C{{\bf C}}    %added

\def\Im{{\mathrm{Im}}}

\def\:{\colon}
\input cyracc.def
\font\tc=wncyr10

\def\i1{\u i}
\newtheorem{theorem}{Theorem}

\title{A Toda lattice in dimension 2 and Nevanlinna theory}
\author{
A.\ Eremenko\thanks{
Supported by NSF grants DMS-0100512 and DMS-0244421.}}
\date{October 20, 2004}
\begin{document}
\maketitle
\begin{abstract} It is shown
how to study the 2-D Toda system for SU(n+1) using
Nevanlinna theory of meromorphic functions
and holomorphic curves. The results generalize recent results
of Jost -- Wang and Chen -- Li.
\end{abstract}

We consider the $2$-dimensional open Toda system for
$SU(n+1)$:
\begin{equation}
\label{toda}
-\frac{1}{2}\Delta u_j=-e^{u_{j-1}}+2e^{u_j}-e^{u_{j+1}},
\quad 1\leq j\leq n,\quad u_0=u_{n+1}=0,
\end{equation}
where $u_j$ are smooth functions in the complex plane $\C$,
and $n$ is a positive integer.

This system was recently studied by several authors
(see the reference list in \cite{jost}). When $n=1$
we obtain the Liouville equation
\begin{equation}
\label{liouville}
-\Delta u_1=4e^{u_1}.
\end{equation}
Jost and Wang \cite{jost} classified all solutions of (\ref{toda})
satisfying the condition
\begin{equation}
\label{summable}
\int_{\C}e^{u_j}<\infty,\quad 1\leq j\leq n.
\end{equation}
In this paper, such classification will be given for a larger
class of solutions, namely those that satisfy
\begin{equation}
\label{growth}
B(r)=\int_{|z|\leq r}e^{u_1}=O(r^K)\quad\mbox{for some}\quad K>0.
\end{equation}
This will also give another proof of the result of Jost and Wang.
Following their suggestion in \cite[p. 278]{jost}, we apply
the Nevanlinna theory of holomorphic curves.

It is well known, and goes back to Liouville \cite{liouville}, 
that the general solution of (\ref{liouville}) has the form
\begin{equation}
\label{spher}
u=u_1=\log\frac{2|f'|^2}{(1+|f|^2)^2},
\end{equation}
where $f:\C\to\CP^1$ is a meromorphic function
with no critical points
(that is $f'(z)\neq 0,\; z\in\C$ and $f$ has no multiple poles).

The general solution of the Toda system can be similarly
described in terms of holomorphic curves $\C\to\CP^n$.
System (\ref{toda}) appears for the first time
in the work on the value
distribution of holomorphic curves
\cite{ahlfors,weyl1,weyl2}. In \cite{calabi},
Calabi proved that every
solution of (\ref{toda}) comes from a holomorphic curve;
this result will be stated precisely below.

Strictly speaking, the present paper contains no new results.
Its purpose is to translate some old results of value distribution
theory to the language of PDE, and thus to bring these results
to the attention of a wider audience. We begin with the simpler
case of the Liouville equation.
\vspace{.1in}

\noindent
{\bf 1. Liouville equation.}
We recall that the order of a meromorphic function $f$
can be defined by the formula
\begin{equation}
\label{order}
\lambda=\limsup_{r\to\infty}\frac{\log A(r,f)}{\log r}.
\end{equation}
where 
$$A(r,f)=\frac{1}{\pi}\int_{|z|\leq r}\frac{|f'|^2}{(1+|f|^2)^2}.$$
If $u$ is related to $f$ by (\ref{spher})
and satisfies (\ref{growth}) then 
\begin{equation}
\label{lamb}\lambda\leq K.
\end{equation}
So (\ref{spher}) establishes a bijective correspondence
between solutions of the Liouville equation satisfying
(\ref{growth}) and meromorphic functions of finite order
without critical points. This class of meromorphic functions
was completely described by F. Nevanlinna in \cite{FNev}.
It coincides with the set of all solutions of
differential equations 
\begin{equation}
\label{schwarz}
\frac{f^{\prime\prime\prime}}{f'}-\frac{3}{2}\left(
\frac{f^{\prime\prime}}{f'}\right)^2=P,
\end{equation}
where $P$ is an arbitrary polynomial.

The expression in the left hand side of this equation is called
the Schwarzian derivative. It is well known (after Schwarz) that
(\ref{schwarz}) is equivalent to the linear differential
equation
\begin{equation}
\label{linear}
w^{\prime\prime}+\frac{1}{2}Pw=0.
\end{equation}
More precisely, solutions $f$ of (\ref{schwarz}) are
ratios $w_1/w_2$ of 
two linearly independent solutions of (\ref{linear}),
and every such ratio is a solution of (\ref{schwarz}). 
This can be verified by a direct computation.

The order $\lambda$ of the meromorphic function $f$ is the same as
the orders of all non-zero solutions of (\ref{linear}),
and it is given by the formula
$$\lambda=(\deg P)/2+1,\quad\mbox{if}\quad P\neq 0.$$
If $P=0$ then all solutions of (\ref{schwarz}) are
fractional-linear functions.

So we have the following recipe for writing all
solutions of (\ref{liouville}) satisfying (\ref{growth}).

\begin{theorem}
Let $u$ be a solution of the Liouville equation
$(\ref{liouville})$ that satisfies $(\ref{growth})$.
Then
$$u=\log \frac{2|f'|^2}{(1+|f|^2)^2},$$
where $f=w_1/w_2$, and the $w_i$ are two linearly independent
solutions of $w^{\prime\prime}+(1/2)Pw=0$, where
$P$ is a polynomial of degree at most $2(K-1)$.

In the opposite direction, if $P$ is a polynomial of degree $d\geq 0$
and $u$ and $f$ are defined as above, then $u$ is a solution
of $(\ref{liouville})$ satisfying $(\ref{growth})$ with $K=d/2+1$.
\end{theorem}

The well known asymptotic behavior of solutions of
the equation (\ref{linear}) permits to derive precise asymptotic
formulas for solutions $u$ of (\ref{liouville}).
We only consider some special cases. First we infer that the order
of growth $\lambda$ of the function $B$ in
(\ref{growth}) can only assume
a discrete sequence of values: $0,\, 1,\, 3/2,\, 2,\, 5/2\ldots.$

If $\lambda=0$, then $P=0$ and $f$ is fractional-linear.
This case was studied by Chen and Li \cite{cl}.
They noticed a curious fact that spherical derivative
of any fractional-linear function always has rotational
symmetry about some center.

If $\lambda=1$,
then $P=\const$, and  $f$ is a fractional-linear function
of $e^{az}$, where $a\neq 0$ is a complex number. So we obtain
a family of elementary solutions of (\ref{toda}).

If $\lambda=3/2$,
then $\deg P=1$, and $u$ can be expressed in terms of the
Airy function.

In the general case, we have the following asymptotic behavior.
Let $P(z)=az^d+\ldots,$ where $d=\deg P$.
The Stokes lines of the equation (\ref{linear}) are defined by
$$\Im (\sqrt{a}z^{d/2+1})=0.$$
They break the complex plane into $d+2$ sectors  of opening
$2\pi/(d+2)$. In each of these sectors, $u(z)$ tends to $-\infty$,
while on the Stokes lines it grows like
$(d/2)\log|z|$.
\vspace{.1in}

\noindent
{\bf 2. Toda system.}
By a holomorphic curve we mean in this paper a holomorphic map
$$f:\C\to\CP^n,$$
where $\CP^n$ is the complex projective space of dimension $n$.
In homogeneous coordinates, a holomorphic curve is
described by a vector-function
$$F:\C\to\C^{n+1},\quad
F=(f_0,f_1,\ldots,f_n)$$ such that $f_j$ are entire
functions without 
common zeros. These entire functions are defined up to a common
multiple which is an entire function without zeros in the plane.
To each holomorphic curve correspond derived curves
defined in homogeneous coordinates as
$$F_k:\C\to\C^{n_k},\quad n_k=\binom{n+1}{k+1},
\quad 1\leq k\leq n,$$
$$F_k=f\wedge f'\wedge\ldots\wedge f^{(k)}.$$
It will be convenient to set
\begin{equation}
\label{bdry}
F_0=F,\quad \mbox{and}\quad F_{-1}=1.
\end{equation}
Notice that $F_n:\C\to\C$ is an entire function; it is equal to the
Wronskian determinant of $f_0,\ldots,f_n$.

$f$ is linearly non-degenerate,
that is its image is not contained
in any hyperplane, if and only if
$F_k\neq 0$ for $1\leq k\leq n$.

The following relations are sometimes called ``local Pl\"ucker
formulas'', see \cite{ahlfors,jost,shabat,weyl2}.
\begin{equation}
\label{plucker}
\Delta\log|F_{k}|^2=4\frac{|F_{k-1}|^2|F_{k+1}|^2}{|F_{k}|^4},\quad
0\leq k\leq n-1.
\end{equation}
Here $|\;|$ is the Euclidean norm.
As $F_n$ is an entire function, we have
\begin{equation}
\label{plucker2}
\Delta\log|F_n|^2=0.
\end{equation}
These relations (\ref{plucker}) and (\ref{plucker2})
hold in the classical sense, that is outside
the zeros of $F_k$.

For a given holomorphic curve $f:\C\to\CP^n$ we set
for $k=1,\ldots,n$:
\begin{equation}
\label{u}
u_k=\log|F_{k-2}|^2-2\log|F_{k-1}|^2+\log|F_{k}|^2+\log 2.
\end{equation}
Notice that the $u_k$ do not change if all $f_k$ are
multiplied by a common factor, so $(u_1,\ldots,u_k)$
depend only on the holomorphic curve $f$ rather than the
choice of its homogeneous representation.

It is  verified by simple computation using (\ref{plucker}),
(\ref{plucker2}) and our conventions (\ref{bdry}) that
these functions $u_k$ satisfy the Toda system (\ref{toda}).

Calabi \cite{calabi} proved the converse statement:
every solution of the system (\ref{toda}) in the plane
arises from a holomorphic curve $f:\C\to\CP^n$,
satisfying 
\begin{equation}
\label{crit}
F_k(z)\neq 0,\; z\in\C,\quad 1\leq k\leq n,
\end{equation}
via (\ref{u}).

Conditions $(\ref{crit})$ are equivalent
to 
\begin{equation}
\label{unram}
F_n(z)\neq 0,\; z\in\C.
\end{equation}

This can be proved by a computation in local coordinates,
see, for example \cite{shabat}. The following arguments are
taken from Petrenko's book \cite{petr}, see also \cite{E2}.

Holomorphic curves that satisfy (\ref{unram}) are called
unramified. 
\vspace{.1in}

\noindent
{\bf Proposition 1} {\em The class of unramified holomorphic
curves $f:\C\to\CP^n$ coincides with the class of curves
that have homogeneous representations of the form
$F=(f_0,\ldots,f_n)$, where $f_1,\ldots,f_n$ is a basis
of solutions of a differential equation
\begin{equation}
\label{equa}
w^{(n+1)}+P_{n}w^{(n)}+\ldots+P_0w=0,
\end{equation}
where $P_j$ are arbitrary entire functions.}
\vspace{.1in}

{\em Proof}. If $f$ is unramified we can choose a homogeneous
representation where $F_n=W(f_0,\ldots,f_n)\equiv 1.$
Then 
$$\left|\begin{array}{cccc} w^{(n+1)}& w^{(n)}&\ldots&w\\
                        f_0^{(n+1)}&f_0^{(n)}&\ldots&f_0\\
                        f_1^{(n+1)}&f_1^{(n)}&\ldots&f_1\\
                        \ldots&\ldots&\ldots&\ldots\\
                        f_n^{(n+1)}&f_n^{(n)}&\ldots&f_n\end{array}
\right|=0$$
is the required equation.

In the opposite direction, suppose that $(f_0,\ldots,f_n)$ is
a fundamental system of solutions of an equation (\ref{equa}).
If the Wronskian $W=W(f_0,\ldots,f_n)$ has a zero,
$W(z_0)=0$, then the columns of the matrix of $W(z_0)$
are linearly dependent, and we obtain a non-trivial
linear combination
$g$ of $f_0,\ldots,f_n$ such that $g^{(k)}(z_0)=0$ for
$k=0,\ldots,n$. As $g$ is a solution
of the same equation (\ref{equa}),
we conclude from the uniqueness theorem
that $g\equiv 0$, but then $W\equiv 0$ which is impossible
because $f_0,\ldots,f_n$ are linearly independent.
This proves the Proposition.
\vspace{.2in}

So we obtained a parametrization of all
smooth solutions of the Toda system in the plane: every
solution has the form (\ref{u}) where $F_k$ are the derived curves
of a curve $f$ whose coordinates form a basis of solutions
of a linear differential equation with entire
coefficients.

Now we turn our attention to condition (\ref{growth}).
We recall that the order of a holomorphic curve can be defined
by equation (\ref{order}) where
$$A(r,f)=\frac{1}{2\pi}\int_{|z|\leq r}\Delta\log|F_0|=
\frac{1}{\pi}\int_{|z|\leq r}\frac{|F_1|^2}{|F_0|^4}.$$
The geometric meaning of $A(r,f)$
is the normalized area of the disc $|z|\leq r$
with respect to the pull-back
of the Fubini-Study area form via $f$.

Let $f$ be a holomorphic curve associated with a solution of
the Toda system via (\ref{u}). Then $B(r)=2\pi A(r,f)$,
and condition (\ref{growth})
implies that $f$ is of finite order, $\lambda\leq K$.
It is known that a holomorphic curve $f$ of finite order
has a homogeneous representation $F=(f_0,\ldots,f_n)$
where the $f_j$ are entire functions of 
finite order, in fact maximum of their orders equals
the order of
$f$ (see, for example \cite{E}).

Let us fix such a representation $(f_0,\ldots,f_n)$. By 
Proposition 1, these entire functions constitute a basis
of solutions of a differential equation
\begin{equation}
\label{diff}
w^{(n+1)}+P_{n}w^{(n)}+\ldots+P_0w=0,
\end{equation}
Now we use a theorem of M. Frei \cite{frei}:
\vspace{.1in}

\noindent
{\bf Proposition 2} {\em If a linear differential equation
of the form $(\ref{diff})$ has $n+1$ linearly independent
solutions of finite order then all $P_j,\; j=0\ldots,n$ are
polynomials.}
\vspace{.1in}

The converse is also true: all solutions of the equation
(\ref{diff}) with polynomial coefficients have finite order.

There is a simple algorithm which permits to find 
orders of solutions of (\ref{diff}).
We follow \cite{wittich}.

Suppose that $P_j\neq 0$ for
at least one $j\in[0,n]$.
Plot in the plane the points with coordinates $(k,\deg P_k-k)$,
for $0\leq k\leq n+1$, $\deg P_{n+1}=0$.
Let $C$ be the convex hull of these
points, and $\Gamma$ the part of the boundary $\partial C$,
which is visible from above. This polygonal line $\Gamma$ is called
the Newton diagram of the equation (\ref{diff}). Then the orders
of solutions are among the
negavive slopes of segments of $\Gamma$. Let
$-\lambda$ be the negative slope of a segment of $\Gamma$,
which has the largest absolute value among all negative slopes
of segments of $\Gamma$. According to a result of
P\"oschl \cite{poschl}, a solution of exact order $\lambda$
always exists. Then $\lambda$ is the order of the 
holomorphic curve defined by (\ref{diff}),
and we have
\begin{equation}
\label{orderl}
\lambda=\max_{0\leq k\leq n}\frac{n+1-k+\deg P_k}{n+1-k}.
\end{equation}
Now we state our final result.
\begin{theorem} Every solution $u$ of the Toda system
in $\C$
that satisfies $(\ref{growth})$ is of the from
$(\ref{u})$, where $F_k$ are the derived curves of a holomorphic
curve whose homogeneous coordinates form a basis of solutions
of
the equation $(\ref{diff})$ with polynomial coefficients.
The degrees of these coefficients, when substituted to
$(\ref{orderl})$, give
$\lambda\leq K$, where $K$ is the constant from $(\ref{growth})$.

Every basis of solutions of any equation $(\ref{diff})$
defines a solution $(u_k)$ of the Toda system by the above rule. 
This solution satisfies $(\ref{growth})$ with $K=\lambda$,
where $\lambda$ is defined by $(\ref{orderl})$. 
\end{theorem}

Some special cases are:

1. If $P_0=\ldots=P_n=0$ in (\ref{diff}),
then a basis of solutions is $(1,z,\ldots,z^n)$
which is the rational normal curve.
Thus we recover the main result of the paper \cite{jost}.

2. Suppose that all $P_j$ are constants. Then solutions of
(\ref{diff}) are generalized exponential sums, and we obtain a class
of explicit solutions of (\ref{toda}).

The asymptotic behavior of solutions of (\ref{toda}) is 
more complicated in the general case $n>1$ than for $n=1$.

The author thanks Misha Sodin who introduced him to
the subject of Toda systems, supplied with the references
and made valuable comments.

{\em Purdue University

150 N University street

West Lafayette IN 47907-2067

USA

eremenko@math.purdue.edu}
\end{document}